\documentclass{article} 

\usepackage{latexsym}
\usepackage{amsfonts}
\usepackage{amssymb}
\usepackage{amsmath}
\usepackage{amscd}
\usepackage{eucal}
\usepackage[all, knot]{xy}
\xyoption{arc}
\xyoption{rotate}
\usepackage{amsthm}
\usepackage{pict2e}
\usepackage{pstricks}
\usepackage{pst-node}
\usepackage{verbatim} 
\usepackage{graphics}
\usepackage{graphicx}
\usepackage{color}   
\usepackage{framed}

\usepackage{bbm}

\usepackage{marvosym}  
\usepackage{stmaryrd}
\usepackage{textcomp}
\usepackage{wasysym}

\psset{unit=1mm,linewidth=0.5pt,
arrowlength=1.1, arrowinset=0.7,arrowsize=4.5pt,
doublesep=1pt,
labelsep=2pt,nodesep=2pt}
  
\def\1{\ensuremath{\mathbbm{1}}}%
\def\2{\ensuremath{\mathbbm{2}}}%

\newcommand{\E}{\ensuremath{\bb{E}}}

\newcommand{\cat}[1]{\ensuremath{\textrm{\bfseries {\upshape {#1}}}}}

\newcommand{\bb}[1]{\ensuremath{\mathbb {#1}}}
\newcommand{\ed}{\end{document}}

\newcommand{\demph}[1]{{\bfseries #1}}

\newcommand{\uscore}{\ensuremath{\underline{\hspace{0.6em}}}}

\newcommand{\tra}{\psset{unit=0.1cm,nodesep=0pt} \pspicture(8,0)
\pcline{->}(1,1.1)(7,1.1) \endpspicture}

\newcommand{\tmapsto}{\psset{unit=0.1cm,nodesep=0pt} \pspicture(8,0) 
\pcline{|->}(1,1.2)(7,1.2) \endpspicture}

\newcommand{\tramap}[1]{\psset{unit=0.1cm,nodesep=0pt,labelsep=2pt} \pspicture(8,4)
\pcline{->}(1,1.1)(7,1.1)\naput{\ensuremath{\scriptstyle{#1}}} \endpspicture}

\newcommand{\tmap}{\tramap}

\newcommand{\ltramap}[1]{\psset{unit=0.1cm,nodesep=0pt} \pspicture(17,4)
\pcline{->}(2.5,1.1)(14.5,1.1)\naput{\ensuremath{\scriptstyle{#1}}} \endpspicture}

\newcommand{\ltmap}{\ltramap}

\newcommand{\dend}{\ensuremath{\hfill \maltese}\end{mydefinition}}

\newcommand{\igc}[1]{\begin{tabular}[c]{c}\includegraphics[width={#1}]}
\newcommand{\igt}[1]{\begin{tabular}[t]{c}\\[-24pt] \includegraphics[width={#1}]}
\newcommand{\ei}{\end{tabular}}

\newtheorem{theorem}{Theorem}[section]

\newtheorem{proposition}[theorem]{Proposition}

\newtheorem{lemma}[theorem]{Lemma}

\newtheorem{definition}[theorem]{Definition}

\newtheorem{example}[theorem]{Example}
\newtheorem{note}[theorem]{Note}
\newtheorem{nonexample}[theorem]{Non-example}
\newtheorem{examples}[theorem]{Examples}
\newtheorem{remarks}[theorem]{Remarks}
\newtheorem{exercise}[theorem]{Exercise}
\newtheorem{question}[theorem]{Question}
\newtheorem{questions}[theorem]{Questions}
\newtheorem{algorithm}[theorem]{Algorithm}
\newtheorem{method}[theorem]{Method}

\newenvironment{mydef}{\begin{definition}\upshape} {\end{definition}}

\newenvironment{myremark}{\begin{remark}\upshape} {\end{remark}}

{ \end{sf}\end{framed}\end{minipage}
\end{center}}

\newtheoremstyle{example}{\topsep}{\topsep}%
     {}
     {}
     {\bfseries}
     {.}
     {8pt}
     {\thmname{#1}\thmnumber{ #2}\thmnote{ #3}}

   \theoremstyle{example}

\newtheorem{remark}[theorem]{Remark}

\newenvironment{prf}{\vspace{1ex}\begin{sloppypar}{\noindent\upshape
{\bfseries Proof.}}} {{\hspace*{\fill}
$\Box$}\end{sloppypar}\vspace{2ex}}

\newcommand{\glob}{
\xy
(-5,0)*+{.}="1";
(5,0)*+{.}="2";
{\ar@/^1pc/^{} "1";"2"};
{\ar@/_1pc/_{} "1";"2"};
{\ar@{=>}^{} (0,2)*{};(0,-2)*{}} ;
\endxy}

\newcommand{\psinv}[6]{\xy
(#1,0)*+{#3}="x";
(#2,0)*+{#5}="y";
{\ar@<.7ex>^{#4} "x"; "y"};
{\ar@<.7ex>^{#6} "y"; "x"};
\endxy
}

\psset{unit=0.1cm,labelsep=1pt,nodesep=3pt}

\begin{document}

\title{A direct proof that the category of 3-computads is not cartesian closed}

\author{Eugenia Cheng\\Department of Mathematics, University of Sheffield \\E-mail: e.cheng@sheffield.ac.uk}

\maketitle

\begin{abstract}
We prove by counterexample that the category of 3-computads is not cartesian closed, a result originally proved by Makkai and Zawadowski.  We give a 3-computad $B$ and show that the functor $\uscore \times B$ does not have a right adjoint, by giving a coequaliser that is not preserved by it.  
\end{abstract}


\setcounter{tocdepth}{2}
\tableofcontents


\section*{Introduction}

Makkai and Zawadowski proved in \cite{mz1} that the category of (strict) 3-computads is not cartesian closed and hence is not a presheaf category.  The result can be considered surprising---for example, the opposite was erroneously claimed in \cite{cj1} (and corrected after Makkai and Zawadowski, in \cite{cj2}).  

The reason is related to the Eckmann-Hilton argument, but the proof given in \cite{mz1}, while having this reason at its heart, uses some sophisticated technology to bring this ``reason'' to fruition---some technical results of \cite{cj1} for Artin glueing, which in turn rely on some technical results of Day \cite{day1}.

In this paper we give a direct counterexample, that is, we give a 3-computad $B$ and a coequaliser
\[
\psset{labelsep=2pt}
\pspicture(32,7)
\rput(0,3){\rnode{A}{$E$}}
\rput(16,3){\rnode{B}{$A$}}
\rput(32,3){\rnode{C}{$C$}}
\psset{nodesep=3pt,arrows=->}
\ncline[offset=3pt]{A}{B}\naput{$$}
\ncline[offset=-3pt]{A}{B}\nbput{$$}
\ncline{B}{C}\naput{$$}
\endpspicture
\]
that is not preserved by the functor $\uscore \times B$, hence $\uscore \times B$ does not have a right adjoint.

The idea behind this counterexample is the same as the idea behind the proof in \cite{mz1}, and the result is, evidently, not new.  However, we believe it is of value to provide this direct argument.  

The root of the problem is that 2-cells having 1-cell identities as source and target do not behave ``geometrically''---by an Eckmann-Hilton argument, horizontal and vertical composition for such cells must be the same and commutative.  Intuitively, this means that cells do not have well-defined ``shape''; a little more precisely, this means for example that if we have 2-cells $a$ and $b$ with identity source and target, then a 3-cell with source $ab \ (=ba)$ cannot have well-defined faces, as we cannot put the putative faces $a$ and $b$ in any order.

This argument obviously does not constitute a proof, but it is the idea at the root of the argument in \cite{mz1} and at the root of the argument we give here.  We begin in Section~\ref{one} by recalling the basic definitions; in Section~\ref{two} we give the counterexample, and in Section~\ref{three} we give the justification.  Experts will only need to read Section~\ref{two}.

Note that unless otherwise stated, all $n$-categories are strict.

\subsection*{Acknowledgements}

I would like to thank and Fran\c cois M\'etayer for asking me for this counterexample.  I would also like to thank Albert Burroni and Yves Guiraud for lively discussions around the subject.

\section{Basic definitions}\label{one}

We begin by recalling the definition of the category of 3-computads.  However, we will only need a small fragment of it for our counterexample, so we will focus on that part.  2-computads are defined by Street in \cite{str7}; the higher-dimensional generalisation is given by Burroni under the name ``polygraphs'' in \cite{bur2} (see also \cite{bat4}).  

The idea is that a 3-computad is a 3-category that is ``level-wise free''.  From another point of view it is the underlying data for a 3-category in which $k$-cells are allowed to have source and target that are pasting diagrams of $(k-1)$-cells, rather than the single $(k-1)$-cells that are the only allowed source and target for globular sets.  Crucially for us, this means in particular that the source and target can be degenerate, that is, identities.

The definition proceeds inductively.  At each dimension we must specify the $k$-cells and then generate pasting diagrams freely in order to specify the boundaries of cells at the next dimension.  This is done using a free 3-category functor and is the technically tricky part of the definition.   However, we will not actually need the full construction of this functor.

\begin{mydef}
A \demph{3-computad} $A$ is given by, for each  $0 \leq k \leq 3$  

\begin{itemize}
\item a set $A_k$ of $k$-cells, and
\item a boundary map $A_k \tra PA_{k-1}$.
\end{itemize}

Here $PA_{k-1}$ denotes the set of parallel pairs of formal composites of $(k-1)$-cells of $A$.  A \demph{morphism of 3-computads} $A \tra B$ is given by, for each $0 \leq k \leq 3$ a morphism
\[f_k: A_k \tra B_k\]
making the obvious squares commute.  We write \cat{3Comp} for the category of 3-computads and their morphisms.

\end{mydef}

In general it is quite complicated to make $P$ precise, but each of the computads involved in our counterexample will have only one 0-cell and no 1-cells.  In this case, the free 2-category on the 2-dimensional data is simply the free commutative monoid on $A_2$ (regarded as a doubly degenerate 2-category).  We use the following terminology.

\begin{mydef}
A 3-computad $A$ is called \demph{2-degenerate} if $A_0$ is terminal and $A_1$ is empty.  Thus by the Eckmann-Hilton argument it consists of
\begin{itemize}
 \item sets $A_2$ and $A_3$, equipped with
\item source and target maps
\[\psset{labelsep=2pt,nodesep=2pt}
\pspicture(15,5)
\rput(0,3){\rnode{a1}{$A_3$}}
\rput(15,3){\rnode{a2}{$A_2^*$}}
\psset{nodesep=3pt,arrows=->}
\ncline[offset=3pt]{a1}{a2}\naput{{\scriptsize $s$}}
\ncline[offset=-3pt]{a1}{a2}\nbput{{\scriptsize $t$}}
\endpspicture
\]
where $A_2^*$ denotes the free commutative monoid on $A_2$.
\end{itemize}

A morphism $A \tra B$ of such 3-computads is given by morphisms
\[A_2 \tmap{f_2} B_2\]
\[A_3 \tmap{f_3} B_3\]
such that the following diagram commutes serially.

\[\psset{unit=0.1cm,labelsep=0pt,nodesep=3pt}
\pspicture(20,20)

\rput(20,18){\rnode{a1}{$A_2^*$}} 
\rput(20,0){\rnode{a0}{$B_2^*$}} 

\psset{nodesep=3pt,arrows=->}
\ncline[offset=0pt]{a1}{a0}\naput[npos=0.5,labelsep=2pt]{\scriptsize{$f_2^*$}} 

\rput(0,0){\rnode{b0}{$B_3$}} 
\rput(0,18){\rnode{b1}{$A_3$}} 
\ncline[offset=0pt]{b1}{b0}\nbput[npos=0.5,labelsep=2pt]{\scriptsize{$f_3$}} 

\ncline[offset=3pt]{b0}{a0}\naput[npos=0.5,labelsep=2pt]{\scriptsize{$s$}} 
\ncline[offset=3pt]{b1}{a1}\naput[npos=0.5,labelsep=2pt]{\scriptsize{$s$}} 

\ncline[offset=-3pt]{b0}{a0}\nbput[npos=0.5,labelsep=2pt]{\scriptsize{$t$}}
\ncline[offset=-3pt]{b1}{a1}\nbput[npos=0.5,labelsep=2pt]{\scriptsize{$t$}}

\endpspicture\]

\end{mydef}

\section{The counterexample}\label{two}

All the 3-computads involved here will be 2-degenerate.  When we check universal properties we will of course need to check them against all computads \emph{a priori}, but we quickly see that the diagrams will ensure 2-degeneracy of any 3-computads involved.

We will write 2-cells as $a,b, \ldots$ and the commutative composition as
\[a.b = b.a.\]
In all that follows, every 3-cell will have a single 2-cell as target, but this is largely to ease the notation; a ``smaller'' counterexample would be possible with empty targets.  

To show that \cat{3Comp} is not cartesian closed we need to show that there exists $B \in \cat{3Comp}$ such that $\uscore \times B$ does not have a right adjoint, so it suffices for $\uscore \times B$ not to preserve all colimits.  So we exhibit a coequaliser
\[
\psset{labelsep=2pt}
\pspicture(36,7)
\rput(0,3){\rnode{A}{$E$}}
\rput(18,3){\rnode{B}{$A$}}
\rput(36,3){\rnode{C}{$C$}}
\psset{nodesep=3pt,arrows=->}
\ncline[offset=3pt]{A}{B}\naput{\scriptsize{$\alpha_1$}}
\ncline[offset=-3pt]{A}{B}\nbput{\scriptsize{$\alpha_2$}}
\ncline{B}{C}\naput{\scriptsize{$\beta$}}
\endpspicture
\]
and a computad $B$ such that the functor $\uscore \times B$ does not preserve it.

\subsection*{Step 1: the coequaliser}

\begin{enumerate}
 \item Let $A$ be the 2-degenerate 3-computad with 2-cells $a_1, a_2, a_3$ and a single 3-cell
\[a_1.a_2 \tmap{f} a_3.\]

\item Let $E$ be the 2-degenerate 3-computad with 2-cells $x,y$ and no 3-cells.

\item Define the morphism $\alpha_1$ by
\[\begin{array}{ccc}
   x & \tmapsto & a_1 \\
   y & \tmapsto & a_3
  \end{array}\]
and define $\alpha_2$ by
\[\begin{array}{ccc}
   x & \tmapsto & a_2 \\
   y & \tmapsto & a_3
  \end{array}\]

\item Thus the coequaliser $C$ simply identifies $a_1$ and $a_2$; it has 2-cells $\bar{a}, a_3$ and a single 3-cell
\[\bar{a}.\bar{a} \tmap{\bar{f}} a_3.\]

\end{enumerate}

\subsection*{Step 2: the functor $\uscore \times B$}

\begin{enumerate}

\setcounter{enumi}{4}

 \item Let $B$ be the 2-degenerate 3-computad (isomorphic to $A$) with 2-cells $b_1, b_2, b_3$ and a single 3-cell
\[b_1.b_2 \tmap{g} b_3.\]

\item $E\times B$ has 2-cells $(x,b_j)$ and $(y,b_j)$ for $j=1,2,3$.  It has no 3-cells.

\item $A\times B$ is the key structure.  It has 2-cells $(a_i, b_j)$ for $i,j = 1,2,3$ and \emph{two} 3-cells
\[\begin{array}{ccc}
   (a_1, b_1). (a_2, b_2) & \ltmap{(f,g)_1} & (a_3, b_3) \\
   (a_2, b_1). (a_1, b_2) & \ltmap{(f,g)_2} & (a_3, b_3) 
  \end{array}\]
This is probably the most interesting part of the argument; we give the full proof later.

\item $C \times B$ has 2-cells $(\bar{a}, b_j)$ for $j = 1,2,3$ and a single 3-cells
\[(\bar{a}, b_1). (\bar{a}, b_2) \ltmap{(\bar{f}, g)} (a_3, b_3).\]

\end{enumerate}

\subsection*{Step 3: non-preservation}

\begin{enumerate}

\setcounter{enumi}{8}

 \item We now examine the coequaliser

\[
\psset{labelsep=2pt}
\pspicture(45,7)
\rput(0,3){\rnode{A}{$E\times B$}}
\rput(25,3){\rnode{B}{$A\times B$}}
\rput(45,3){\rnode{C}{$P$}}
\psset{nodesep=3pt,arrows=->}
\ncline[offset=3pt]{A}{B}\naput{\scriptsize{$\alpha_1 \times 1$}}
\ncline[offset=-3pt]{A}{B}\nbput{\scriptsize{$\alpha_2 \times 1$}}
\ncline{B}{C}\naput{$$}
\endpspicture
\]
and show that it is not isomorphic to $C \times B$.

Now the morphism $\alpha_1 \times 1$ is given by
\[\begin{array}{ccc}
   (x,b_j) & \tmapsto & (a_1,b_j) \\
   (y,b_j) & \tmapsto & (a_3,b_j)
  \end{array}\]
and $\alpha_2 \times 1$ by
\[\begin{array}{ccc}
   (x,b_j) & \tmapsto & (a_2,b_j) \\
   (y,b_j) & \tmapsto & (a_3,b_j)
  \end{array}\]
Thus the coequaliser $P$ simply identifies $(a_1,b_j)$ with $(a_2,b_j)$ for each $j$.  So it has 2-cells which we may call $(\bar{a}, b_j)$ and $(a_3, b_j)$ (which is to be expected as the coequaliser is preserved up to 2 dimensions).

$P$ has two distinct 3-cells
\[\begin{array}{ccc}
   (\bar{a}, b_1). (\bar{a}, b_2) & \ltmap{(f,g)_1} & (a_3, b_3) \\
   (\bar{a}, b_1). (\bar{a}, b_2) & \ltmap{(f,g)_2} & (a_3, b_3). 
  \end{array}\]
Since $C \times B$ has only one 3-cell it is clear that $C \times B$ is not isomorphic to this coequaliser $P$, that is, $\uscore \times B$ does not preserve the original coequaliser.

Note that the canonical factorisation
\[P \tra C \times B\]
identifies the 3-cells $(f,g)_1$ and $(f,g)_2$.

\end{enumerate}

\section{Universal properties}\label{three}

In this section we check all the universal properties required for the counterexample.  In principle we only need to check the 3-cells, as 2-computads form a presheaf category so we know that the lower dimensions behave pointwise.  However we include the full argument for completeness, and because it is straightforward.

\begin{lemma}
 The product $A \times B$ is as given in the previous section, with the obvious projections.
\end{lemma}

\begin{prf}
 We exhibit its universal property.  Consider a 3-computad $Y$ and morphisms
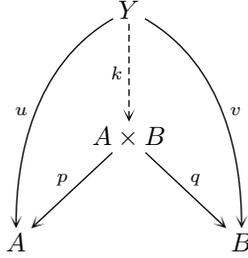
\[
\psset{unit=0.1cm,labelsep=2pt,nodesep=2pt}
\pspicture(0,40)(40,72)


\rput(5,40){\rnode{a}{$A$}}  
\rput(20,54){\rnode{ab}{$A \times B$}}  
\rput(20,71){\rnode{c}{$Y$}}  

\rput(35,40){\rnode{b}{$B$}}  

\ncline{->}{ab}{a} \nbput{{\scriptsize $p$}}
\ncline{->}{ab}{b} \naput{{\scriptsize $q$}}

\nccurve[angleA=210,angleB=90,ncurvA=0.8,ncurvB=0.6]{->}{c}{a}\nbput{{\scriptsize $u$}}

\nccurve[angleA=-30,angleB=90,ncurvA=0.8,ncurvB=0.6]{->}{c}{b}\naput{{\scriptsize $v$}}

\psset{linestyle=dashed, dash=2.5pt 1.5pt}
\ncline{->}{c}{ab} \nbput{{\scriptsize $k$}}

\endpspicture
\]

We seek to exhibit a unique factorisation $k$ as shown. On 0-, 1- and 2-cells, $A \times B$ is just a product, so we define the factorisation at these dimensions as for products ie 
\[k(t) = (u(t), v(t)).\]
Note in particular that $A$ and $B$ have no 1-cells, so for the morphisms $u$ and/or $v$ to exist, $Y$ cannot have any 1-cells either.  So this map respects boundaries trivially.

%
%
%
%
%
%
%
%
%
%
%

We now discuss the factorisation on 3-cells.  Let $e$ be a 3-cell in $Y$.  Now $A$ and $B$ have only one 3-cell each, $f$ and $g$ respectively.  So we must have
\[\begin{array}{ccc}
   u(e) &=& f \\
   v(e) &=& g 
  \end{array}\]
thus e must have boundary as follows
\[y_1. y_2 \tmap{e} y_3\]
for some 2-cells $y_1, y_2, y_3 \in Y$.  Then since the action of $u$ and $v$ respect the boundary of $e$ we know $y_3$ must be sent to $a_3$ and $b_3$ respectively.  However considering the source there is some ambiguity as the product is commutative, so for each of $u$ and $v$ there are two possibilities---either the subscripts are left the same, or they are switched.  That is, on ordered pairs the action of $u$ is
\[\begin{array}{rccc}
   \mbox{either} & (y_1, y_2) & \tmapsto & (a_1, a_2) \\
\mbox {or} & (y_1, y_2) & \tmapsto & (a_2, a_1)
  \end{array}\]
and similarly the action of $v$ is
\[\begin{array}{rccc}
   \mbox{either} & (y_1, y_2) & \tmapsto & (b_1, b_2) \\
\mbox {or} & (y_1, y_2) & \tmapsto & (b_2, b_1).
  \end{array}\]
There are thus 4 cases, but in each case $k(e)$ is uniquely determined to be either $(f,g)_1$ or $(f,g)_2$ by the condition that $k$ preserves boundary.  Explicitly, $k(e)$ is specified by examining the action of $u$ and $v$ as shown by the following table.
\[\begin{array}{cc|ccc}
&&&v&\\
&& (y_1, y_2) \tmapsto (b_1, b_2) && (y_1, y_2) \tmapsto (b_2, b_1)\\[4pt]
\hline &&&& \\[-4pt]
& (y_1, y_2) \tmapsto (a_1, a_2) & (f,g)_1 && (f,g)_2 \\[-4pt]
u&&& \\[-4pt]
& (y_1, y_2) \tmapsto (a_2, a_1) & (f,g)_2 && (f,g)_1
\end{array}
\]

%
%
%
%
%
%
%
%
%
%
%
%
%
%

\end{prf}

The other products follow similarly, but more easily. It remains to check the universal properties of the two coequalisers in question, which is much more straightforward.

Consider a diagram 

\[
\psset{labelsep=2pt}
\pspicture(40,27)
\rput(0,23){\rnode{A}{$E$}}
\rput(20,23){\rnode{B}{$A$}}
\rput(40,23){\rnode{C}{$C$}}
\psset{nodesep=3pt,arrows=->}
\ncline[offset=3pt]{A}{B}\naput{\scriptsize{$\alpha_1$}}
\ncline[offset=-3pt]{A}{B}\nbput{\scriptsize{$\alpha_2$}}
\ncline{B}{C}\naput{\scriptsize{$q$}}

\rput(40,6){\rnode{Y}{$Y$}}
\ncline{B}{Y}\nbput{\scriptsize{$u$}}
\ncline[linestyle=dashed,dash=3pt 2pt]{C}{Y}\naput{\scriptsize{$k$}}

\endpspicture
\]
with $u\alpha_1 = u\alpha_2$.  We seek a unique factorisation $k$ as shown.  

\begin{itemize}

\item On 0-cells: $A$ and $C$ only have one 0-cell each; writing each as $\ast$ we must have $k(\ast) = u(\ast) \in Y$.

\item On 1-cells: $A$ and $C$ have no 1-cells, so as before $Y$ cannot have any either.  

\item On 2-cells: To make the triangle commute we must put
\[\begin{array}{ccl}
   k(\bar{a}) &=& u(a_1) \ [= u(a_2)]\\
  k(a_3) &=& u(a_3).
  \end{array}\]
This respects boundaries as all 2-cells involved are degenerate.

\item On 3-cells: To make the triangle commute, we must have $k(\bar{f}) = u(f)$.  This respects boundaries, by our definition of $k$ on 2-cells.

\end{itemize}

The other coequaliser proceeds in the same way, but with two 3-cells.

\begin{myremark}
Note that this sort of counterexample cannot arise for 2-computads, as 2 is the lowest dimension of cell for which the Eckmann-Hilton argument can be used.  Note also that this problem does not arise for weak 3-computads as weak identity 1-cells impede the Eckmann-Hilton argument on degenerate 2-cells.  This difference between the commutativity of degenerate 3-cells in weak and strict structures also arises in \cite{cm1}.
\end{myremark}


\ed